\documentclass[12pt]{article}

\usepackage{times}

\usepackage{graphicx}
\usepackage[final]{pdfpages}
\newenvironment{sciabstract}{%
\begin{quote} \bf}
{\end{quote}}

\newcounter{lastnote}

\title{Periodically fighting shake, rattle and roll}

\author
{M. Brun$^{1}$, A.B. Movchan$^{2\ast}$, I.S. Jones$^{3}$\\
\\
\normalsize{$^{1}$ Department of Structural Engineering, University of Cagliari,
Cagliari I-09123, Italy}\\
\normalsize{$^{2}$ Department of Mathematical Sciences, University of Liverpool, } \\
\normalsize{ Peach Street,
Liverpool, L69 3BX, UK}\\
\normalsize{$^{3}$ School of Engineering,
     John Moores University, Liverpool L3 3AF, UK }\\
\\
\normalsize{$^\ast$To whom correspondence should be addressed; E-mail:  abm@liv.ac.uk}
}

\date{}

\begin{document}


\baselineskip24pt


\maketitle


\begin{sciabstract}

How easy is it to suppress shake, rattle and roll in a long bridge or a skyscraper?
Most practical  structures are designed so that long wave resonance vibrations can be avoided.
However, there are recent examples, such as the Millennium Bridge in London or the Volga Bridge in Volgograd,
which show that unexpected external forces may result in large scale unwanted shake and rattle.
Full scale alteration of a bridge (or a skyscraper) would not be considered as an acceptable option, unless the structure has collapsed.

Can we fix this by examining a representative part of the structure only and making small lightweight changes?
We will do it here and illustrate an idea linking the engineering analysis to elastic waveguides.

%This approach is generic and may also be applied to the mechanics of earthquake-resistant civil engineering structures.

\end{sciabstract}

%\section*{Introduction}
\paragraph*{  Background of the problem.}  %Introduction}

Figs. 1A and 1B show the troubled automobile bridge across the river Volga in Volgograd.
The large deformation of the upper deck was completely unexpected.
As described in the Daily Mail  (see \cite{DailyMail}) and in IMechanica professional discussion
forum \cite{imech8280}, large vibrations have been invoked by relatively small external forces.
The bridge was inaugurated in October 2009, and in May 2010 a long wave resonance vibration caused sections of the
bridge to bend as  shown in Figs. 1B and 1C.
A fundamental alteration of the bridge on its foundation appears to be problematic and no efficient technical
solutions have been proposed to date.
The model illustrated here gives a solution of this challenging problem and, being generic, also leads to
design methodologies  for vibrating slender structures such as skyscrapers and earthquake-resistant systems.
The analysis complements the models of lateral vibration induced by ``crowd-synchronization'' and
``balancing pedestrians'' for the Millennium Bridge discussed in \cite{StrAbrMcREckOtt05,San08,McD09}.

Although the general theory of waves in periodic structures is a classical topic described
in many textbooks (see, for e\-xam\-ple, \cite{brillouin53} and \cite{kittel53}),
%{\color{red}  until recently  its development for elastic systems was limited. The bulk of the work  was}
% {\color{blue}$\rightarrow$ } {\color{blue} 
 until recently  the bulk of the work  was %}
related to problems of electromagnetism and acoustics
as confirmed by the extensive bibliography  on photonic band gap structures
\cite{PhononicReferences} compiled in 2008 by J.P. Dowling.

%%%%%%%%%%%%%%%%%%%%%%%%%%%%%%%%%%%%%%%%%%%%%%%%%%%%%%%%%%%%%%%%%%%%%
\begin{figure}[!h]
\begin{centering}
\includegraphics[width=16cm]{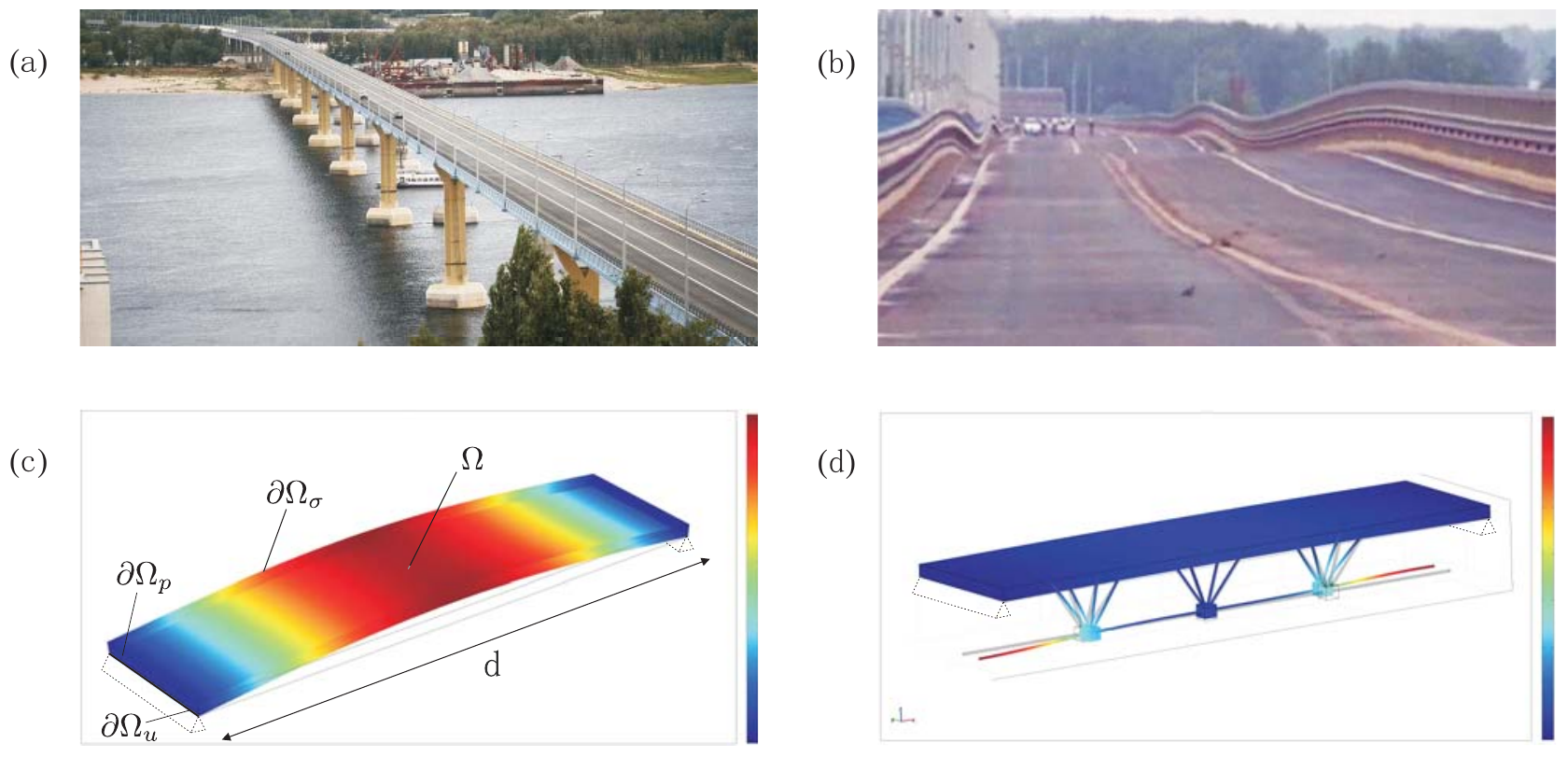}
\vspace{0.2cm} \caption{Volga bridge, ((A), (B) HTB Volgograd News).
Flexural deformation of the main upper body of the bridge: (B) real structure, (C) simplified numerical model,
(D) modified 3D structured wave guide with lightweight resonators.}
\label{fig1}
\end{centering}
\end{figure}
%%%%%%%%%%%%%%%%%%%%%%%%%%%%%%%%%%%%%%%%%%%%%%%%%%%%%%%%%%%%%%%%%%%%%

In the recent years, the theory of elastic %Bloch-Floquet
waves in periodic systems has
%{\color{red}  received substantial attention,} 
%{\color{blue}$\rightarrow$ } 
%{\color{blue}  
received increasing attention, %}
with applications ranging from the design of
elastic filters and polarisers to the modelling of earthquake-resistant structures.
Analytical models for phononic crystals, having  periodic
structures, have been developed in
\cite{MovNicMcP97,PouMovMcPNicAnt00,ZalMovPouMcP02,MovGueMovMcP07,MovMovGueMcP07}.
Transmission problems for arrays of elastic structured stacks were studied in \cite{PlaMovMcPMov02,PlaMovMcP03,PlaMovMcPMov03}.
These incorporate the comparative analysis of the filtering properties of %Bloch-Floquet
elastic waves in doubly periodic media and the transmission properties for the corresponding singly periodic stack structure.
Propagation and dispersion of waves in plate structures have been analysed in \cite{MovMovMcP07,McPMovMov09,MovMcPMovPou09,JonMovGei10}; in particular,
Green's functions and localised vibration modes within plate structures are included in \cite{McPMovMov09}.
Heterogeneous lattices were used in
\cite{BruGueMovBig10,BruMovMov10} to build models of structured interfaces, which possess specially designed transmission properties.
These models also included the means of controlling  dispersion properties.
Effects of disorder on propagation of waves in structured media as well as localisation have been addressed in \cite{PlaMovMcPMov03,McPMovMov09}, {%\color{red} 
where the multipole Rayleigh method and a novel recurrence algorithm were used to obtain the transmission characteristics of a composite stack incorporating multiple arrays of circular inclusions/voids.}

{ %\color{blue}  I HAVE MOVED THIS PARAGRAPH TO THIS POSITION.
Smart composite structures have been designed and built for a range of applications in optics and problems of
electromagnetism to encompass phenomena such as negative refraction and the cloaking of finite sized objects.
Complex models have been constructed in \cite{Pen99,Smi05,PenSchSmi05}.
The theoretical findings have also been verified experimentally \cite{SheSmiSch01,SmiPenWil04,SchMocJusCumPenStaSmi04,ErgSteBrePenWeg10}.
Specially designed composite structures, referred to as metamaterials, may channel  waves around finite sized obstacles (cloaking) or create the effect of negative refraction which can be further used in focussing of electromagnetic waves by flat interfaces.
Analysis of the geometrical transforms used in modelling of such materials is presented in \cite{NicZol09}.
Modelling of elastic waves brings new challenging questions in the design of metamaterials. % structures.
In particular, a geometric transform was used in \cite{BruGueMov08} to model an ``invisibility cloak'' for elastic waves.}

In this paper, we consider a problem where elastic waves have to be channeled around some part of an engineering structure,
for example the main deck of a bridge.  { %\color{red} 
The background  work of the authors, which  addresses analysis of waves interacting with structured interfaces, is included in \cite{BruGueMov08, JonMovGei10, BruGueMovBig10}.}
Although, the problem is very different from the formulations of electromagnetism, there is a direct link with the design of
metamaterial structures in \cite{ErgSteBrePenWeg10}.
Physically, the interpretation of the model can be thought of in terms of a ``by-pass" for the
%damaging
{ %\color{red} 
undesired} elastic waves.
The method of design is non-trivial, { %\color{red} 
and the analysis presented here is original.}
%Rather than {\color{red} studying} a full-scale finite structure and computing its numerous modes,
%we reduce the model to the analysis of elastic waves in a periodic waveguide.
%This
{ %\color{red} 
Our} approach provides accurate estimates { %\color{red} 
for frequencies of  standing waves within the  elastic structure} and the design for the by-pass of elastic waves, so that shake, rattle and roll is no more.

%\section*{Periodic waveguide}
\paragraph*{{ %\color{red} 
Physical model for a} periodic waveguide.}

Instead of looking into resonance modes of a finite slender body such as a bridge or a tall building, we consider a periodic waveguide
and then analyse dispersion properties of { %\color{red} 
Bloch-Floquet} elastic waves.
A typical span of the bridge between two neighboring pillars is shown in Fig. 1C and will be used as an elementary cell of the periodic structure.

%For convenience, we shall
{  %\color{red} 
{\it The geometrical notations} are introduced as follows.} We denote by $\Omega$  the domain occupied by the elementary cell {  %\color{red} 
of the periodic waveguide.}
The boundary of $\Omega$ is divided into the constrained part $\partial \Omega_u$, which is in contact
with the supporting pillars, the traction free part $\partial \Omega_\sigma$ and the contact region between consecutive unit cells $\partial \Omega_p$. {  %\color{red} 
These notations are shown in Fig. 1C, Fig. 2A and Fig. 6H.}

The elastic displacement is assumed to be time-harmonic of angular frequency $\omega$,
and its amplitude ${\bf u}$ satisfies the Lam\'e  equation
\begin{equation}
\mu \Delta {\bf u} +(\lambda+\mu)\nabla \nabla\cdot{\bf u}+\rho\omega^2{\bf u}=0 \quad\mbox{in~} \Omega,  ~%$$%
\label{eq1a}
\end{equation}
where $\mu$ and $\lambda$ are the Lam\'e {%\color{red} 
elastic moduli} %constants
and $\rho$ { %\color{red} 
is the mass density.} The boundary conditions are
\begin{equation}
{\bf u} = 0, \quad \mbox{on}\quad \partial \Omega_u, \label{bc1}
\end{equation}
{%\color{red} 
representing the fixed part of the boundary,}
and { %\color{red} 
on the remaining boundary subjected to the traction boundary conditions we set}
\begin{equation}
{\bf t}^{(n)}({\bf u}) = 0, \quad \mbox{on}\quad \partial \Omega_\sigma, \label{bc2}
\end{equation}
where ${\bf t}^{(n)}({\bf u})$ is the vector of tractions.
%Any propagating
{%\color{red} 
A Bloch-Floquet wave}, characterised by the wave number $k$,  would also satisfy  the quasi-periodicity condition within the elementary cell
\begin{equation}
{\bf u}({\bf x}+d{\bf e}^{(1)})={\bf u}({\bf x})e^{ikd},~ \label{bfc}
\end{equation}
where $d$ is the period of the structure.
Subject to modification of boundary conditions, a similar formulation can be used to describe a time-harmonic vibration of a tall building.
Analysis of dispersion properties of elastic waves characterizes the
%energy propagation
{%\color{red} 
group velocity} for different values of the
frequency $\omega$.
In the low frequency range, this analysis also provides information about stop bands,
defined as intervals of frequencies where no waves can propagate through the structured wave\-guide.
One would also learn about standing waves, which may exist within the slender engineering structure and,  in turn, one would be able to suppress undesired vibrations or channel waves around some parts of the structure such as the upper deck of a long bridge.

%\section*{Two-dimensional illustrative example}
\paragraph*{Two-dimensional   { %\color{red} 
implementation.}}    %illustrative example}

Firstly, we provide an illustration of a periodic waveguide for a two-dimensional analogue of the model {%\color{blue} 
of the upper deck} shown in Fig. 1C, which is given as a periodically constrained {%\color{red} 
slender elastic solid, that can be approximated in the low frequency range as an elastic beam.
 We refer to \cite{KozMazMov}
%Kozlov {\it et al. } (1999)
for technical details of the asymptotic analysis, which leads to the lower-dimensional {\it beam type} approximations of solutions of the problem (\ref{eq1a})--(\ref{bfc}) for the Lam\'e system in the desired range of frequencies.}
Within the elementary cell, the beam is simply supported and such a system  possesses low frequency flexural modes.
We show with a simple representative example (see Fig. 2), that such flexural modes can be suppressed by the addition of a lightweight periodic system of resonators.
The analysis of dispersion properties of elastic waves enables one to tune the system accurately and to choose the correct design.

%%%%%%%%%%%%%%%%%%%%%%%%%%%%%%%%%%%%%%%%%%%%%%%%%%%%%%%%%%%%%%%%%%%%%
\begin{figure}[!h]
\begin{centering}
\includegraphics[width=11cm]{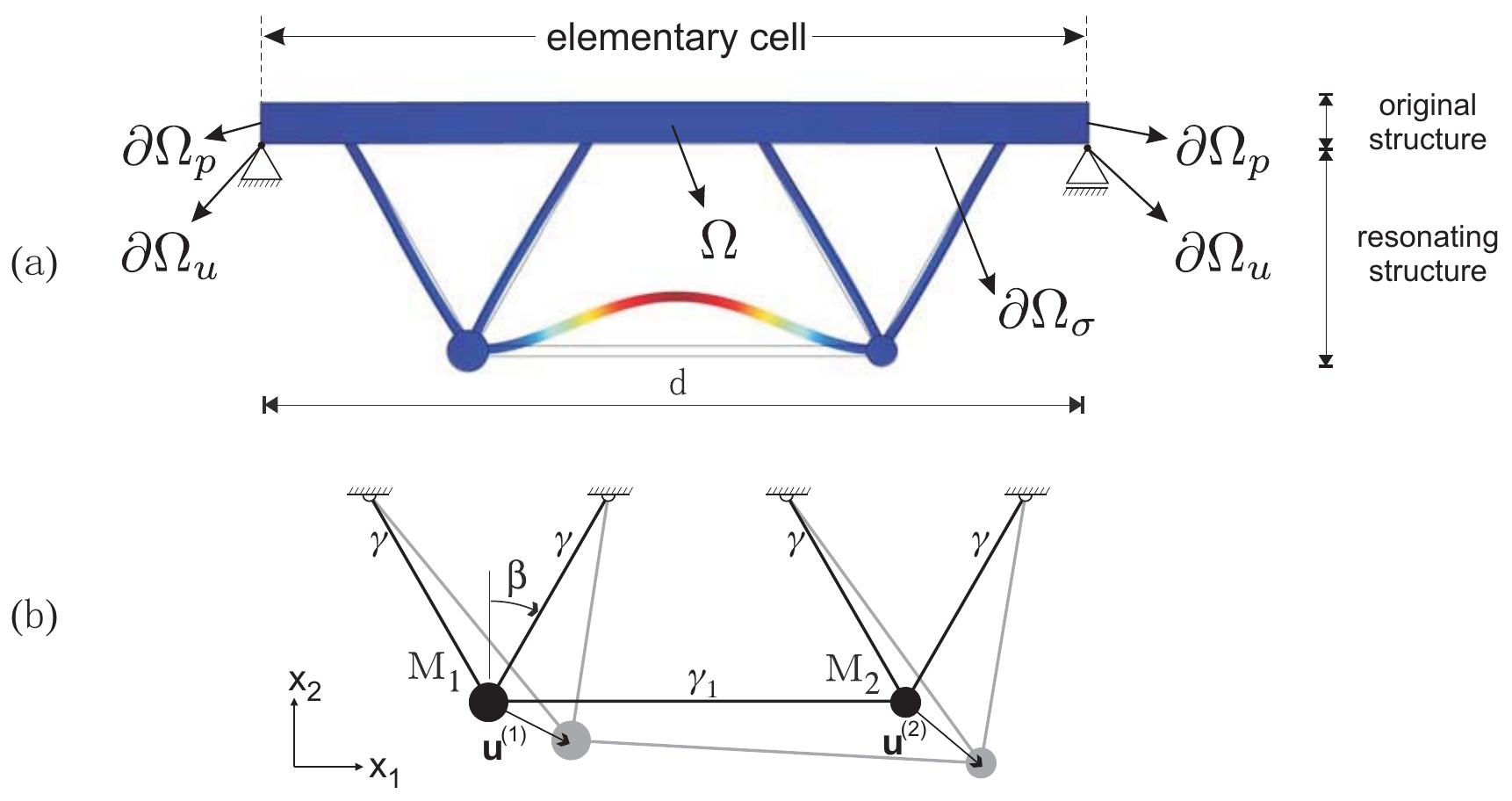}\\
\vspace{0.2cm} \caption{Two-dimensional elementary cell with the system of resonators.
(A) The first eigenmode.
(B) Truss model used for the analytical estimates { %\color{blue} 
of the eigenfrequencies}: each mass $M_i$ undergoes the horizontal and vertical displacements $u^{(i)}_1$ and $u^{(i)}_2$,
respectively, with $i=1,2$. Parameter values are:
$d=4$ m; the thickness $s=0.2$ m; the radii of the disks are $0.1$ m and $0.075$ m; $h_1=2$ m, $h_2=1$ m, $\beta=\pi/6$;
the longitudinal stiffness coefficients are $\gamma=0.14$ GPa, $\gamma_1 = 0.018$ GPa; the main plate has mass density $\rho = 7850$ kg/m$^3$
and shear modulus $\mu=80$ GPa; the disks and the elastic links have mass density $\rho_M=7850$ kg/m$^3$ and  $\rho_\gamma = 200$ kg/m$^3$,
respectively.}
\label{fig2}
\end{centering}
\end{figure}
%%%%%%%%%%%%%%%%%%%%%%%%%%%%%%%%%%%%%%%%%%%%%%%%%%%%%%%%%%%%%%%%%%%%%

For flexural vibrations of the periodically constrained beam, the dispersion diagram,
presenting the normalised frequency ${\cal F}=fd/v$ as a function of $k d$, is shown in Fig. 3A (Here $f = \omega/2 \pi$ is the frequency, $d$ is the period, $v$ is the speed of the shear wave in the upper deck of the bridge).

%The corresponding eigenmode at ${\cal F}=0.0358$,
%representing a flexural wave, is similar to the one shown in Fig. 1C.

{%\color{red} 
The flexural mode, that is similar to the one shown in Fig. 1C,  corresponds to ${\cal F}=0.0358$. Our aim is to use the data from the dispersion diagram to tune the elastic system so that the above flexural mode is suppressed.  }

%%%%%%%%%%%%%%%%%%%%%%%%%%%%%%%%%%%%%%%%%%%%%%%%%%%%%%%%%%%%%%%%%%%%%
\begin{figure}[!h]
\begin{centering}
\includegraphics[width=14cm]{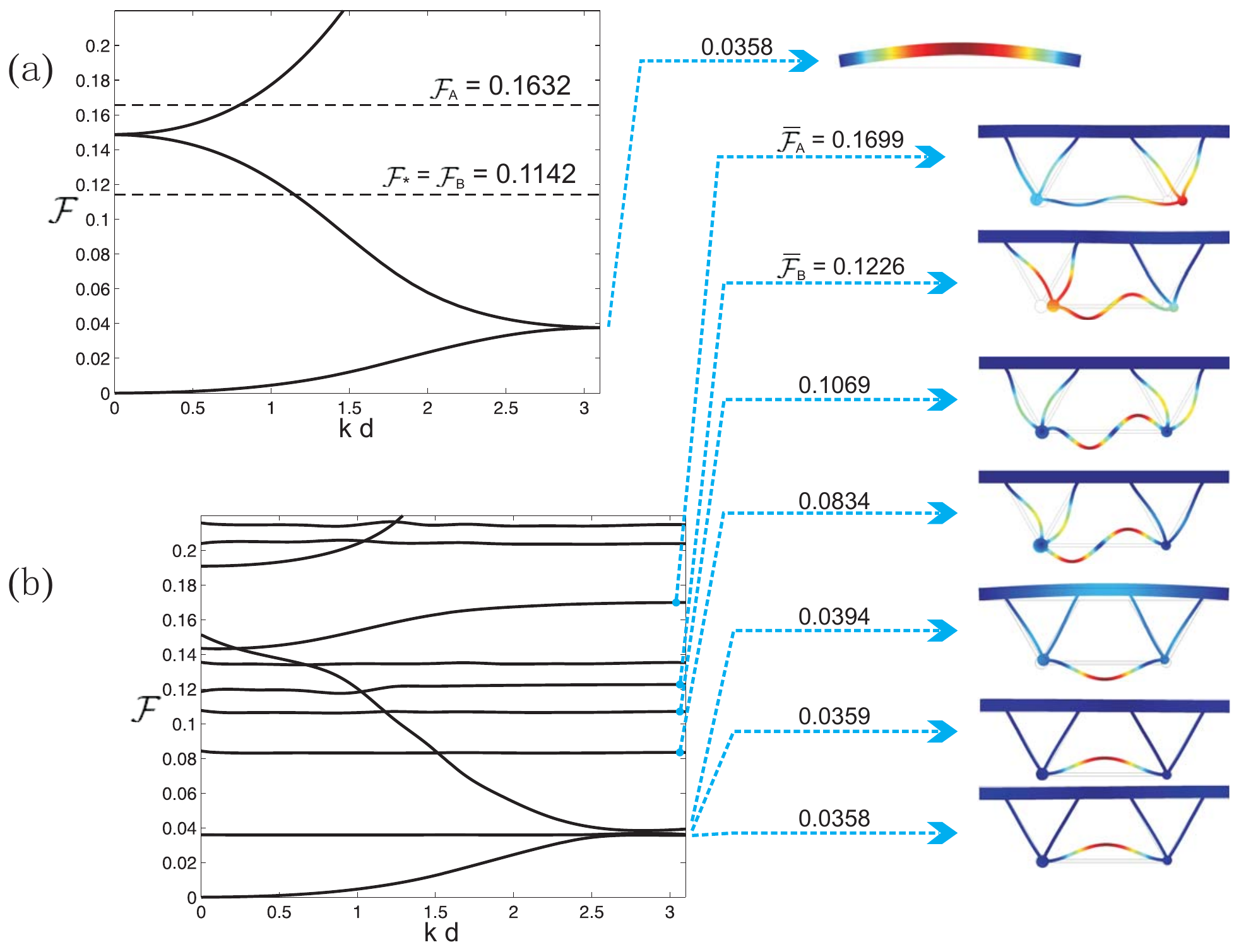}\\
\vspace{0.2cm} \caption{Dispersion diagrams and corresponding eigenmodes for the 2D structured wave guide (finite element computations).
(A) The original beam structure and the eigenmode to suppress. ${\cal F}_A$ and ${\cal F}_B$ correspond to the analytical values given
in eqn. (\ref{eigenfreq})a.
(B) The modified structure and the first several eigenmodes with corresponding normalized frequencies.
$\overline{\cal F}_A$ and $\overline{\cal F}_B$ are finite element results for the analytical frequencies ${\cal F}_A$ and ${\cal F}_B$.
}
\label{fig3}
\end{centering}
\end{figure}
%%%%%%%%%%%%%%%%%%%%%%%%%%%%%%%%%%%%%%%%%%%%%%%%%%%%%%%%%%%%%%%%%%%%%

Now, we introduce a lightweight modification of the original structure which would %eliminate
alter the long-wave flexural oscillations. Such a modification does not involve any major alteration
of the structure such as a change of the boundary conditions in the contact region between the deck of
the bridge and the supporting pillars or variation of the bending stiffness of the upper deck.

For the two-dimensional configuration, the proposed design is shown in Fig. 2, where a connected system
of lightweight resonators is attached to the main body of the bridge.
The overall mass of a single resonator is much smaller compared to the total mass of the elementary cell of the structure.
Also, this system  has a relatively low stiffness ($\gamma$, $\gamma_1$) of  additional elastic links compared to the stiffness of the
main deck of the bridge.
In turn,  the eigenfrequencies of such a re\-so\-na\-tor system
can be placed sufficiently close to the low eigenfrequencies of the original macrostructure.

The new set of dispersion curves is shown in Fig. 3B, which contains stop bands,
i.e. the frequency intervals where no propagating elastic waves exist.
We note a narrow band gap around the normalised frequency  ${\cal F}=0.0359$ and a wider band gap at normalised frequencies around ${\cal F}=0.1699$.

{%\color{red} 
It is vital to pay attention to  the vibration modes near the boundaries of the band gaps mentioned and emphasise  that they  do not involve large flexural vibrations of the main body of the bridge, as required in the proposed design.}

It is important to be able to tune the system correctly to match the band gap frequency with the frequency of an unwanted vibration of the upper deck of the bridge. This task can be achieved by using the {%\color{red} 
analytical} approximations described below.

%\subsection*{Tuning of the  elastic system} %eigenfrequencies of the energy localizers}
\paragraph*{%\color{red} 
Analytical tuning of the vibrating elastic system.} %eigenfrequencies of the energy localizers}

Analytical estimates of the eigenfrequencies for a class of standing waves provide a %preliminary
{%\color{red} 
convenient} design {%\color{red} 
tool} for the resonating
structure and are accompanied by a numerical finite element simulation for an elementary cell comprising a computation of
the eigenmodes and providing a dispersion diagram for elastic waves within the structured wave\-guide.

The approximation scheme is illustrated
%on the simple
{%\color{red} 
for the configuration displayed}
%example
 in Fig. 2B, where the system of resonators is represented
as a low-frequency truss structure with two masses $M_1$ and $M_2$ and
linear elastic links with longitudinal stiffness values $\gamma$ for
the diagonal connections and $\gamma_1$ for the horizontal one.
Assuming that the upper deck is fixed, the low eigenfrequencies of vibration involving
translational motion of masses $M_1$ and $M_2$ are given in the simple analytical form
\begin{eqnarray}
\nonumber
& f_{A,B}^2=\frac{1}{8 \pi^2}\left(\frac{1}{M_1}+\frac{1}{M_2}\right)
\left[\gamma_1+2\gamma\sin^2\beta\pm
\sqrt{\gamma_1^2+4\left(\frac{M_1-M_2}{M_1+M_2}\right)^2\gamma\sin^2\beta(\gamma_1+\gamma\sin^2\beta)}
\right] &,\\
&f_C^2=\frac{\gamma}{2\,M_1 \pi^2}\cos^2\!\beta, ~\quad
f_D^2=\frac{\gamma}{2\, M_2 \pi^2}\cos^2\!\beta,&
\label{eigenfreq}
%\nonumber
\end{eqnarray}
The frequencies $f_{A, B}$ $(f_B<f_A)$ correspond to  horizontal translational motions of the interconnected masses $M_1$ and $M_2$,
whereas the frequencies $f_{C}, f_{D}$ describe the vertical motion of the masses $M_1$ and $M_2$, respectively.
We use the notation $f_*=\min \{f_B, f_C, f_D\}$, representing the quantity
which serves as an upper bound for a cluster of frequencies corresponding to standing waves associated with
rotational motion of the resonators and flexural vibration of the inertial links.
Some of these standing waves are shown to correspond to flat bands on the dispersion diagram in Fig. 3B.
For our particular choice of material parameters, {%\color{red} 
we have} $f_* = f_B$. We also note that
$f_A$ gives an analytical estimate for the position of the second band gap displayed in Fig. 3B.

A supplementary  finite element computation for an elementary cell of the structure, addressed in Fig. 3B,
allows for a full interaction between the elastic upper deck of the bridge and the built-in resonator structure,
which acts as the e\-ner\-gy localizer and hence reduces significantly the amplitude of vibration of the upper
deck of the bridge, as shown in Fig. 2A.

Comparative analysis of the dispersion diagrams in Figs. 3A and 3B shows formation of band gaps and a cluster
of standing waves for the structure containing lightweight resonators.
It appears that some of the  waves of sufficiently low group velocity, illustrated in Fig. 3B, correspond to
vibrations of low frequency near the boundaries of the band gaps.
The change of the inertial and stiffness properties of the high-contrast lightweight resonator structure
would alter the position of the stop bands and hence can be used as the means of control of propagating waves
of other frequencies within the periodic elastic system.

Further fine tuning of the system may involve, for example, a very slight change of stiffness $\gamma_1$ of the horizontal bar.
This would lead to an accurate match of a band gap frequency with  the frequency of vibration of the upper deck of the unaltered bridge.
In the proposed design, vibrations of the upper deck of the bridge become negligibly small compared to the vibration of the additional resonator structure, as shown in Fig. 2A, and this has been achieved without any alterations of junction conditions between the bridge deck and supporting
pillars or any change in the stiffness of the main deck.
A full three-dimensional simulation is discussed in the text below.

%\section*{Three-dimensional simulation}
\paragraph*{Three-dimensional simulation.}

We use the ideas described above to give %a possible
{%\color{red} 
the required } design solution for the
``Volga bridge problem", discussed in \cite{DailyMail, imech8280}.
In the three-dimensional computational model, we consider two cases: the original bridge structure, which
possesses a standing flexural wave as shown in Fig. 1C, and the new design including a system of lightweight
resonators illustrated in Fig. 1D, {%\color{red} 
where the flexural wave in question has been suppressed.}

As above, instead of referring to a full scale finite length bridge,
we consider a waveguide, with emphasis on the dispersion properties
of elastic waves within a periodic structure.

%%%%%%%%%%%%%%%%%%%%%%%%%%%%%%%%%%%%%%%%%%%%%%%%%%%%%%%%%%%%%%%%%%%%%
\begin{figure}[!h]
\begin{centering}
\includegraphics[width=14cm]{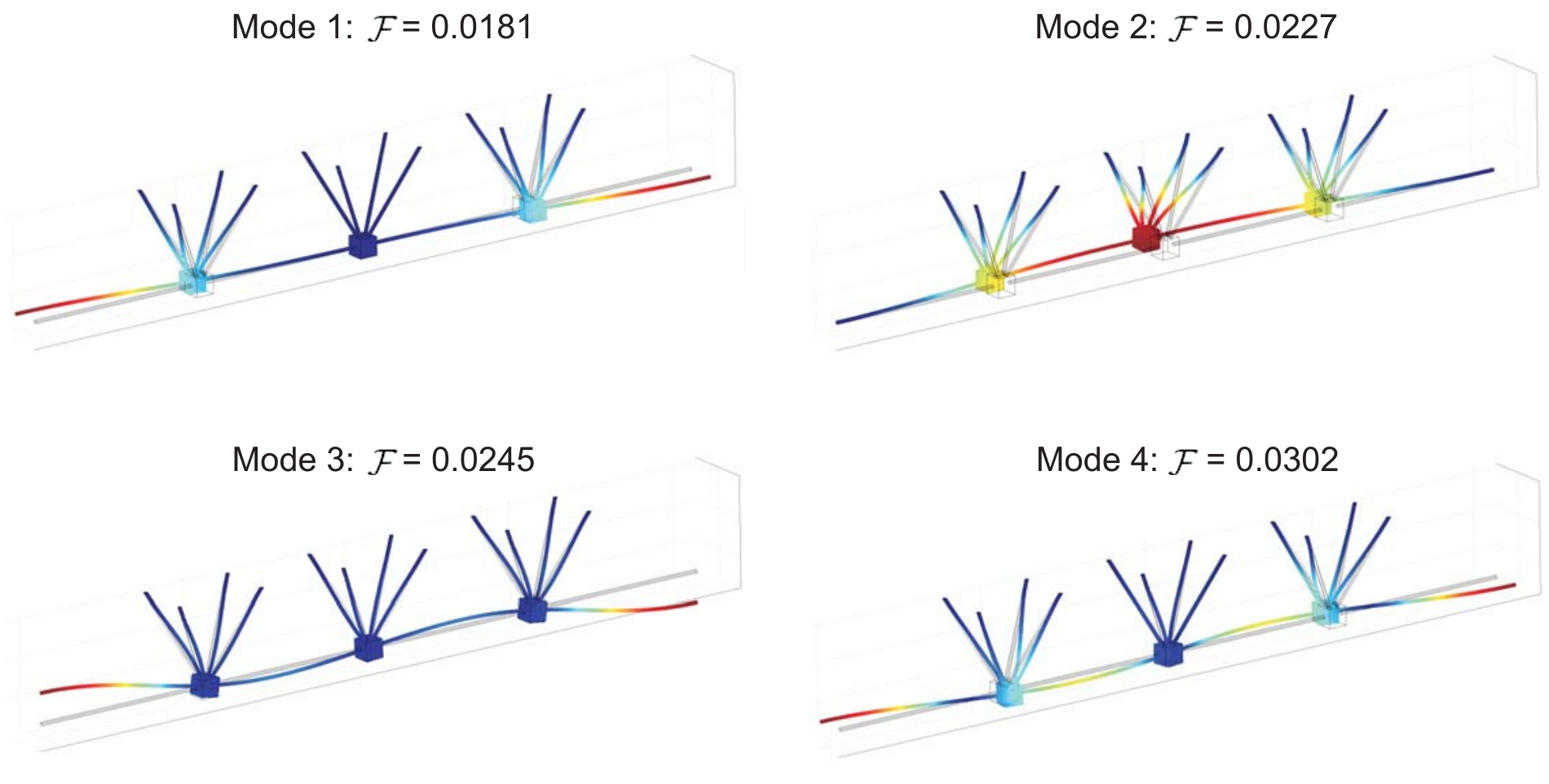}\\
\vspace{0.2cm} \caption{First four vibration modes of the resonant structure.
The first mode, at normalised frequency ${\cal F}= fd/v= 0.0181$, is designed to eliminate the flexural
vibration of the deck of the bridge.
In the computation $d= 4$ m and $v=\sqrt{\mu/\rho}=3194$ m/s is the shear wave speed in the upper deck. }
\label{fig4}
\end{centering}
\end{figure}
%%%%%%%%%%%%%%%%%%%%%%%%%%%%%%%%%%%%%%%%%%%%%%%%%%%%%%%%%%%%%%%%%%%%%

A simple analysis of the resonator structure on its own (as shown in Fig. 4), with the assumption that the surfaces
of contact with the bridge deck are fixed, gives an accurate estimate of frequencies for a class of standing waves
within the modified bridge after installation of the periodic system of resonators.
Fig. 4 also shows the four eigenmodes within the frequency interval $0<{\cal F}<0.0303$.
The physical and stiffness parameters of the resonators are chosen so that one of the
eigenfrequencies of the lightweight resonator structure of Fig. 4 matches the frequency ${\cal F}=0.0182$
of the standing wave of the unaltered bridge.
When the two frequencies in question are sufficiently close to each other
the combined structure will change its dynamic response within the required low frequency range.
Embedding of a periodic system of low-frequency resonators leads  to the formation  of the required cluster of standing waves, where the amplitude of vibration of the bridge deck becomes negligibly small.
The proposed mechanical system can be tuned to filter waves of desired frequencies and hence to
eliminate vibration of the upper deck of the bridge within the specified frequency range.

%%%%%%%%%%%%%%%%%%%%%%%%%%%%%%%%%%%%%%%%%%%%%%%%%%%%%%%%%%%%%%%%%%%%%
\begin{figure}[!h]

\begin{centering}
\includegraphics[width=14cm]{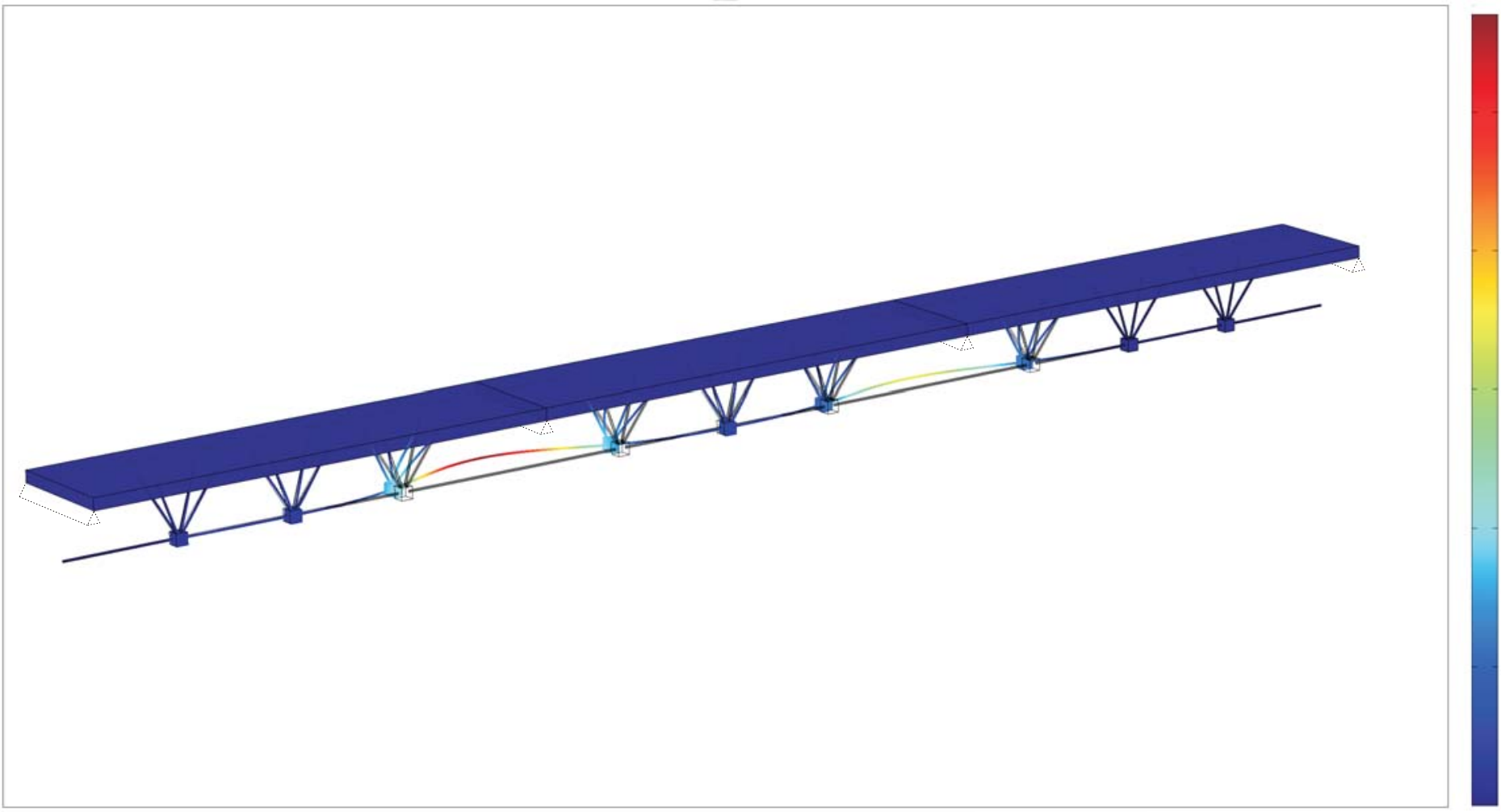}\\
\vspace{0.2cm} \caption{Modified bridge: the flexural vibration of the main body has been suppressed.
The built-in resonators have taken on the vibrational motion.}
\label{fig5}
\end{centering}
\end{figure}
%%%%%%%%%%%%%%%%%%%%%%%%%%%%%%%%%%%%%%%%%%%%%%%%%%%%%%%%%%%%%%%%%%%%%

The computation shows that for the unaltered bridge,
the fundamental flexural mode corresponds to the flexural vibrations of the
upper deck of the bridge at a low normalised frequency ${\cal F}= 0.0182$.
After embedding of the appropriately tuned lightweight structure of the resonators
the eigenfrequency in question has been slightly shifted, and the eigenmodes of vibration have been
dramatically changed.

In Fig. 5 we show a finite structure representing a sufficiently large section of the modified bridge,
where the left and the right boundaries of the deck of the bridge are simply supported.
By analysing the vibration mode of the modified structure, we confirm that the motion of the main body
of the bridge is negligibly small compared to the motion within the system of resonators.
This gives the desired design, which eliminates the flexural vibration of the upper deck
of the bridge without any major alteration of the structure in terms of mass, stiffness or boundary conditions.

%\section*{Suppression of lateral vibrations of a skyscraper }
\paragraph*{Suppression of lateral vibrations of a skyscraper. }

The generic method based on the a\-nal\-y\-sis of waves in a periodic waveguide is now applied to such a system as  a tall building, or `skyscraper', composed of  several copies of an elementary cell, as shown in Fig. 6 panel $(2)$.
{%\color{red} 
The displacement vector is assumed to be time-harmonic, with the radian frequency $\omega$ and the amplitude ${\bf u}$, which satisfies the Lam\'e system (\ref{eq1a}) and the boundary conditions (\ref{bc1}) and (\ref{bc2}) on the   fixed foundation and on the boundary, which is free of tractions, respectively.}

%%%%%%%%%%%%%%%%%%%%%%%%%%%%%%%%%%%%%%%%%%%%%%%%%%%%%%%%%%%%%%%%%%%%%%%%%%%%%%%%%%%%%%%%%%%%%%%%%%%%%%%%%%%%%%%%%%
\begin{figure}[!ht]
\begin{centering}
\includegraphics[height=14.cm]{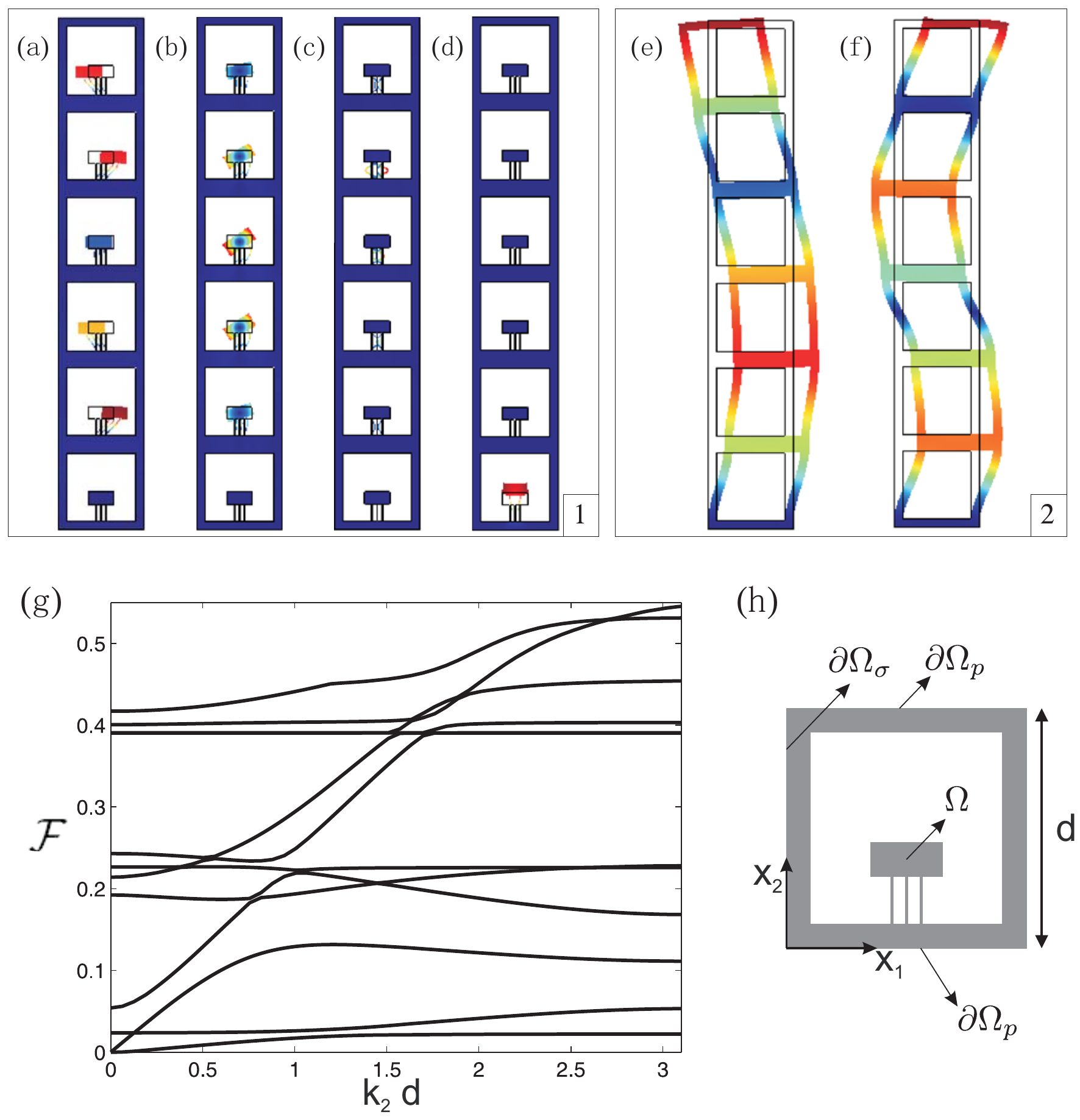}\\
\vspace{0.2cm}
\caption{Eigenmodes of the multi-storey buildings. The horizontal foundation of the structure is fixed,
whereas all the remaining boundary is traction free.
Panel 1: first four eigenmodes of the multi-structure,
each showing energy absorption by the resonators. Normalised eigenfrequencies ${\cal F}=f d/v$:
(A) $0.0225$; (B) $0.2254$; (C) $0.3936$; (D) $0.4446$.
Panel 2: first two eigenmodes of the frame structure. Normalised eigenfrequencies:
(E) $0.0144$; (F) $0.0260$.
(G) The dispersion diagram showing normalised frequency ${\cal F}$ versus $k_2 \, d$.
(H) Square elementary cell with dimension $d=1$ m. The external structure has thickness $s=0.1$ m,
the resonator is composed of a rectangular solid with dimensions $0.3$ m $\times$ $0.14$ m and
three thin ligaments with dimensions $0.2$ m $\times$ $0.01$ m.
Material parameters are: shear modulus $\mu = 80$ GPa, Poisson ratio $\nu = 0.28$, mass density $\rho = 7850$ Kg/m$^3$.}
\label{BuildingModes}
\end{centering}
\end{figure}
%%%%%%%%%%%%%%%%%%%%%%%%%%%%%%%%%%%%%%%%%%%%%%%%%%%%%%%%%%%%%%%%%%%%%%%%%%%%%%%%%%%%%%%%%%%%%%%%%%%%%%%%%%%%%%%%%%

%Here, we
{%\color{red} 
We} introduce a multi-scale resonator {%\color{red} 
within every elementary cell, as shown in Fig. 6H,} and %we
make use of the analysis of elastic waves within a periodic
plain strain structure. % whose elementary cell is shown in Fig. 6H.
The resonator consists of a mass connected to the main frame of the elementary cell by a set of parallel thin
ligaments, whose overall flexural stiffness is sufficiently small.

As a design target, we would like to alter the {%\color{red} 
flexural} vibration modes of the original structure, within a % certain
{%\color{red} 
predefined} range of frequencies,
%providing
{%\color{red} 
to provide} a set of vibration modes which do not involve a large amplitude motion of the walls of the building.
Moreover, the ``thin-legged resonators" may move and hence absorb the energy of the system.
This is clearly a task of high practical importance, which may appear in problems such as the design of earthquake-resistant structures.

The lowest translational and rotational eigenfrequencies of the thin-legged resonators can be estimated
analytically \cite{JonMovGei10,KozMazMov} and the design of the resonators can be chosen
to place these eigenfrequencies within the required interval corresponding to the modes of the original structure
in panel $(2)$ of Fig. 6.

%Now, we
% draw attention to the
{%\color{red} 
The} dispersion diagram of Fig. 6G {%\color{red} 
has been} computed for the periodic system with the unit cell given in Fig. 6H. This includes a group of flat bands corresponding to standing waves related to vibration of resonators within
the periodic structure.
We emphasize the connection between the standing waves associated with vibrations of the thin-legged resonator
and groups or clusters of eigenmodes of the finite structure, some of which are shown in Fig. 6, panel $(1)$.
For such frequencies, the motion of the sides of the skyscraper is negligibly small.
We note that an ad-hoc estimate of such frequencies is difficult for any realistic structure,
and hence the appropriate spectral %Bloch-Floquet
analysis on a single unit cell is desirable, in particular,
in the design process.

The two low frequency modes of vibration of the ``multi-storey'' structure without built-in resonators in Fig. 6 panel $(2)$
show that the whole macro-structure behaves like a homogenized beam.
It is noted that the building, without the system of resonators, will experience lateral vibrations of significant amplitude
and, in practical configurations, this may lead to an overall structural failure.
In Fig. 6 panel $(2)$ we show the structure vibrating at normalised frequencies ${\cal F}=0.0144,  0.0260$.
These values represent two neighboring frequencies at the ends of the interval containing ${\cal F}=0.0225$, which is the first eigenfrequency of the thin-legged resonator, as shown in Fig. 6 panel $(1)$.

The alteration of the vibration modes of panel $(2)$ of Fig. 6 is achieved by introduction of the resonators, as shown in
Fig. 6 panel $(1)$.
This is fully consistent with the  idea of  the modified design for the Volga Bridge discussed earlier.

We pay particular attention to clusters of normalised eigenfrequencies around the values $0.0225$, $0.2254$, $0.3936$ and $0.4446$,
corresponding to standing waves.
For such waves, the exterior boundary of the elementary cell moved with very small amplitude compared to the
amplitude of vibrations of the thin-legged resonators. The four frequencies outlined above correspond to lateral motion of the resonator,
rotational mode of the resonator, localised vibrations of thin legs within the resonator and finally, the longitudinal translational mode of the
resonator, respectively.
The associated eigenmodes are shown in Fig. 6 panel $(1)$
where the color map corresponds to the  total displacement.
In Fig. 6A we have a mode dominated by the lateral translational motion of the resonators;
in Fig. 6B we observe the rotational motion of the resonators; Fig. 6C gives an example from
a large cluster of eigenmodes corresponding to the case of localised vibrations of thin legs
within the resonators.
Finally, in Fig. 6D we observe a localisation near the foundation of the structure (namely
only one resonator located at the ground level is moving, whereas the motion of  the remaining
structure has a negligibly small amplitude).
This last example shows an exponential localisation within a structured waveguide;
in this case, the skyscraper is treated as a waveguide and the embedded structure of resonators
creates a stop band preventing the waves of certain
frequencies from propagating along such a waveguide.
In all of these computations, the amplitude of vibration of the sides
of the large structure is very small compared to that of the
individual resonators.
Hence, if the solid is subjected to an external impact, such resonance modes
can be initiated and %hence
the energy can be absorbed into the resonators and then
dissipated via the desired mechanism, for example damage of the thin legs.
Even though you may not be able to eat your lunch off the resonating tables the building will survive the earthquake.

%\section*{Dissipation of energy and final remarks}
\paragraph*{Dissipation of energy and final remarks}

In the re-designed bridge structure shown in Fig. 5 the waves have been channeled
through the system of resonators, away from the upper deck
of the bridge. For the case of a skyscraper, inclusion of a system of resonators has led
to suppression of low-frequency lateral vibrations.
Although the models discussed above did not incorporate energy dissipation, additional
viscous dampers attached to the resonators can be considered as a feasible development in
a practical implementation.

The main emphasis of the paper is on the analysis of
waves in a periodic system, rather than resonance modes of a finite solid. It is shown that the analysis of the elementary cell
of the periodic waveguide provides insight and required numerical data for the design of the full-scale engineering system.
The ideas presented in this paper are generic, and the bridge problem in addition to the structural
modification of a skyscraper,
can be considered as possible practical implementations.
The proposed methods of analysis would equally apply
to other slender structures in civil engineering and problems of structural design.

%\section*{Acknowledgments}
\paragraph*{Acknowledgments}
We would like to thank Professor Sir John Pendry and Professor Vladimir Shalaev for reading the manuscript and for their support. The authors thank Dr. Anna Matveyeva (Volgograd State University) for having given direct evidence of the structural problem of the bridge
and Professor Paglietti
(University of Cagliari) for the stimulating discussions.
We would like to thank Professor R.C. McPhedran of the University of Sydney for valuable discussion and suggestions on the text of the article.
The financial support of the Visiting Professorship Program 2010 (ABM) of the Regione Autonoma della Sardegna
as well as financial support of Research Centre in Mathematics and Modelling (MB) of the University of
Liverpool are gratefully acknowledged.


\begin{thebibliography}{99}

\bibitem{DailyMail} {\small http://www.dailymail.co.uk/news/worldnews/article-1280919/Russian-bridge-bounces-feet-Volga-River.html}

\bibitem{imech8280} http://imechanica.org/node/8280

% \bibitem{arup} http://www.arup.com/millenniumBridge/
\bibitem{StrAbrMcREckOtt05} S. H. Strogatz, D. M. Abrams, A. McRobie, B. Eckhardt and E. Ott, {\it Nature}, {\bf 438}, 43-44, (2005)

\bibitem{San08} K. Sanderson, {\it Nature}, doi:10.1038/news.2008.1311, (2008).

\bibitem{McD09} J. H. G. McDonald, {\it Proc. R. Soc. A}, {\bf 465}, 1055–1073, (2009).

\bibitem{brillouin53} L. Brillouin, {\it Wave propagation in periodic structures: electric filters and crystal lattices}, (Dover, New York, 1953).

\bibitem{kittel53} C. Kittel, {\it Introduction to solid state physics}, (J. Wiley \& Sons, New York, 1953).

\bibitem{PhononicReferences} The Photonic and Sonic Band Gap Metamaterial Bibliography, 2008 available from: http://phys.lsu.edu/$\sim$jdowling/pbgbib.html


\bibitem{MovNicMcP97} A. B. Movchan, N. A. Nicorovici and R. C. McPhedran, {\it Proc. R. Soc. A}, {\bf 453}, 1958, 643-662, (1997).

\bibitem{PouMovMcPNicAnt00} C. G. Poulton, A. B. Movchan, R. C. McPhedran, N. A. Nicorovici and Y. A. Antipov, {\it Proc. R. Soc. A}, {\bf 456}, 2002, 2543-2559, (2000).

\bibitem{ZalMovPouMcP02} V. V. Zalipaev, A. B. Movchan, C. G. Poulton and R. C. McPhedran, {\it Proc. R. Soc. A}, {\bf 458}, 2026, 2327-2347, (2002).

\bibitem{MovGueMovMcP07} N. V. Movchan, S. Guenneau, A. B. Movchan and R. C. McPhedran, {\it Physica B}, {\bf 394}, 2, 281-284, (2007).

\bibitem{MovMovGueMcP07} A. B. Movchan, N. V. Movchan, S. Guenneau and R. C. McPhedran, {\it Proc. R. Soc. A}, {\bf 463}, 2080, 1045-1067, (2007).

\bibitem{PlaMovMcPMov02} S. B. Platts, N. V. Movchan, R. C. McPhedran and A. B. Movchan, {\it Proc. R. Soc. A}, {\bf 458}, 2024, 1887-1912, (2002).

\bibitem{PlaMovMcP03} S. B. Platts, N. V. Movchan, R. C. McPhedran and A. B. Movchan, {\it J. Eng. Mater. Technol.}, {\bf 125}, 1, 2-6, (2003).

\bibitem{PlaMovMcPMov03} S. B. Platts, N. V. Movchan, R. C. McPhedran and A. B. Movchan, {\it Proc. R. Soc. A}, {\bf 459}, 2029, 221-240, (2003).

\bibitem{MovMovMcP07} A. B. Movchan, N. V. Movchan and R. C. Mcphedran, {\it Proc. R. Soc. A}, {\bf 463}, 2086, 2505-2518, (2007).

\bibitem{McPMovMov09} R. C. McPhedran, A. B. Movchan and N. V. Movchan, {\it Mech. Mat.}, {\bf 41}, 4, 356-363, (2009).

\bibitem{MovMcPMovPou09} N. V. Movchan, R. C. McPhedran, A. B. Movchan and C. G. Poulton, {\it Proc. R. Soc. A}, {\bf 465}, 2111, 3383-3400, (2009).

\bibitem{JonMovGei10} I. S. Jones, A. B. Movchan and M. Gei, {\it Proc. R. Soc. A}, doi: 10.1098/rspa.2010.0319.

\bibitem{Pen99} J. B. Pendry, {\it Science}, {\bf 285}, 5434, 1687-1688, (1999)

\bibitem{Smi05} D. R. Smith, {\it Science}, {\bf 308}, 5721, 502-503, (2005)

\bibitem{PenSchSmi05} J. B. Pendry, D. Schurig and D. R. Smith, {\it Science}, {\bf 312}, 5781, 1780-1782, (2006)

\bibitem{SheSmiSch01} R. A. Shelby, D. R. Smith and S. Schultz, {\it Science}, {\bf 292}, 5514, 77-79, (2001)

\bibitem{SmiPenWil04} D. R. Smith, J. B. Pendry and M. C. K. Wiltshire, {\it Science}, {\bf 305}, 5685, 788-792, (2004)

\bibitem{SchMocJusCumPenStaSmi04} D. Schurig, J. J. Mock, B. J. Justice, S. A. Cummer, J. B. Pendry, A. F. Starr and D. R. Smith {\it Science}, {\bf 314}, 5801, 977-980, (2006)

\bibitem{ErgSteBrePenWeg10} T. Ergin, N. Stenger, P. Brenner, J. B. Pendry and M. Wegener, {\it Science}, {\bf 328}, 5976, 337-339, (2010)

\bibitem{NicZol09} A. Nicolet and F. Zolla, {\it Science}, {\bf 323}, 5910, 46-47, (2009)

\bibitem{BruGueMov08} M. Brun, S. Guenneau and A. B. Movchan, {\it Appl. Phys. Lett.}, {\bf 94}, 6, 061903, (2009).

\bibitem{BruGueMovBig10} M. Brun, S. Guenneau, A. B. Movchan and D. Bigoni, {\it J. Mech. Phys. Solids}, {\bf 58}, 9, 1212-1224, (2010).

\bibitem{BruMovMov10} M. Brun, A. B. Movchan and N. V. Movchan, {\it Cont. Mech. Therm.}, doi:10.1007/s00161-010-0143-z, (2010).

\bibitem{KozMazMov} V. Kozlov, V. Maz'ya, A. Movchan, {\it Asymptotic Analysis of Fields in Multi-structures.} Oxford University Press (1999).


%\bibitem{Images1} {\small http://images.yandex.ru/yandsearch?text=volgograd+bridge}
% \bibitem{youtube1} {\small http://www.youtube.com/watch?v=5smsMzA\_xII}
% \bibitem{youtube2} {\small http://www.youtube.com/watch?v=84sO7UiHcmE}
% \bibitem{youtube3} {\small http://www.youtube.com/watch?v=lgtjG5kUfGs}
% \bibitem{youtube4} {\small http://www.youtube.com/watch?v=yWeHUXPi0Jc}

%\bibitem{youtube5} {\small http://www.youtube.com/watch?v=AQT-t4HQmGo}

% http://www.youtube.com/watch?v=IEl59FCBE-k&playnext_from=TL&videos=W3HjLpPlC0E&feature=grec_index
% http://www.youtube.com/watch?v=MnZrgiQa7WA&feature=watch_response
% http://www.youtube.com/watch?v=WEQrt_w7gN4
% http://www.youtube.com/watch?v=5smsMzA_xII

%Do not include separate BibTeX files; if BibTeX is used, paste the output here.
\end{thebibliography}
\end{document}